\documentclass[11pt]{amsart}
\usepackage{amsxtra}
\usepackage{graphicx}
\usepackage{amssymb,mathabx}
\usepackage{color}
\usepackage[all]{xy}
\usepackage{comment}
\usepackage{nameref}

\usepackage{multicol}
\usepackage{tabularx}

\theoremstyle{plain}

\newcommand{\cleqn}{\setcounter{equation}{0}}
\newcommand{\clth}{\setcounter{theorem}{0}}
\newcommand {\sectionnew}[1]{\section{#1}\cleqn\clth}

\newtheorem{theorem}{Theorem}[section]
\newtheorem{lemma}[theorem]{Lemma}
\newtheorem{definition-theorem}[theorem]{Definition-Theorem}
\newtheorem{proposition}[theorem]{Proposition}
\newtheorem{corollary}[theorem]{Corollary}
\newtheorem{definition}[theorem]{Definition}
\newtheorem{example}[theorem]{Example}

\newtheorem{remark}[theorem]{Remark}
\newtheorem{conjecture}[theorem]{Conjecture}
\newtheorem{notation}[theorem]{Notation}

\newtheorem{question}[theorem]{Question}
\newcommand \bth[1] { \begin{theorem}\label{#1} }
\newcommand \ble[1] { \begin{lemma}\label{#1} }

\newcommand \bpr[1] { \begin{proposition}\label{#1} }
\newcommand \bco[1] { \begin{corollary}\label{#1} }
\newcommand \bde[1] { \begin{definition}\label{#1}\rm }
\newcommand \bex[1] { \begin{example}\label{#1}\rm }
\newcommand \bre[1] { \begin{remark}\label{#1}\rm }
\newcommand \bcj[1] { \begin{conjecture}\label{#1}\rm }

\newcommand \bnota[1] { \begin{notation}\label{#1}\rm }
\renewcommand {\eth} { \end{theorem} }
\newcommand {\ele} { \end{lemma} }

\newcommand {\epr} { \end{proposition} }
\newcommand {\eco} { \end{corollary} }
\newcommand {\ede} { \end{definition} }
\newcommand {\eex} { \end{example} }
\newcommand {\ere} { \end{remark} }
\newcommand {\ecj} { \end{conjecture} }

\newcommand {\enota} { \end{notation} }

\def \Rset {{\mathbb R}}         
\def \Cset {{\mathbb C}}

\def \Zset {{\mathbb Z}}

\def \Nset {{\mathbb N}}

\def \Pset {{\mathbb P}}

\def \cp   {{\Pset^1}}

\def \al {\alpha}

\def \Lie { {\mathrm{Lie \,}} }

\def \g  {\mathfrak{g}}   

\DeclareMathOperator \GCR { {\mathcal{ECR}} }

\DeclareMathOperator \Ker { {\mathrm{Ker}} }

\renewcommand \max { {\mathrm{max}} }

\newcommand{\Rvw}{\mathcal{R}_v^w}
\newcommand{\oRvw}{\mathring{\mathcal{R}}_v^w}

\begin{document}
\title[Powers of projective lines inside $G/B_+$]{Powers of projective lines inside full flag varieties}
\author[Q\"endrim Gashi]{Q\"endrim R. Gashi}
\address{Department of Mathematics, University of Maryland, College Park, MD 20742}
\email{qgashi@umd.edu}
\author[Milen Yakimov]{Milen Yakimov}
\address{
Department of Mathematics, Northeastern University, Boston, MA 02115 and International Center for Mathematical Sciences, Institute of Mathematics and Informatics \\
Bulgarian Academy of Sciences \\ 
Acad. G. Bonchev Str., Bl. 8 \\
Sofia 1113, Bulgaria}
\email{m.yakimov@northeastern.edu}
\date{}
\keywords{Flag varieties, Richardson varieties, Poisson degeneracy loci, toric varieties}
\subjclass[2020]{Primary 14M15; Secondary 53D17, 14L35, 05E18}
\begin{abstract} Up to now, symplectic foliations of Poisson varieties have been primarily studied from a set theoretic perspective by describing their symplectic leaves individually as orbit closures. In the direction of obtaining algebro-geometric characterizations, for the standard Poisson structure $\pi$ on the full flag variety $G/B_+$ of a complex semisimple algebraic group $G$, we prove that 
the irreducible components of the reduced part of the 0-degeneracy locus of $(G/B_+, \pi)$ are precisely the inclusion-maximal $T$-stable complete Poisson submanifolds of $(G/B_+, \pi)$ that are isomorphic to powers of $\cp$. This characterization is based on a theorem that explicitly classifies the closed Richardson varieties of $G/B_+$ isomorphic to $(\cp)^d$. 
In turn, this theorem is based on a general theorem for arbitrary $T$-stable subvarieties of full flag varieties $G/B_+$ isomorphic to $(\cp)^d$.  
\end{abstract}
\maketitle
\sectionnew{Introduction}
\label{sec1}
For a connected simply connected complex semisimple algebraic group $G$, 
with a pair of opposite Borel subgroups $B_\pm$ and corresponding maximal torus $T:= B_+ \cap B_-$, its \emph{open Richardson varieties} are the intersections 
\begin{equation}
\label{oRvw}
\oRvw := B_{-}vB_{+}/B_{+} \cap B_{+}wB_{+}/B_{+}
\end{equation}
of opposite Schubert cells in the full flag variety $G/B_+$
for pairs of Weyl group elements $v \leq  w \in W$. The \emph{closed Richardson varieties} $\Rvw$ are their Zariski closures 
\[
\Rvw := \overline{\oRvw} = 
\bigsqcup_{v \leq v' \leq w' \leq w} \mathring{\mathcal{R}}_{v'}^{w'}.
\]
They play a key role in the theory of Kazhdan--Lusztig polynomials \cite{Deodhar}, standard monomial theory \cite{BrionLakshmibai}, total 
positivity \cite{MR}, Poisson geometry \cite{GY}, Postnikov's positroid cells \cite{KLS}, mirror symmetry \cite{RW}, and many other settings.
Cluster structures on open Richardson varieties were conjectured to exist in \cite{Leclerc}, and this was proved in \cite{CGGSS,GLSB,SSB}. 

In recent years, there has been a strong interest in investigating Richardson varieties from the perspectives of toric geometry and Poisson geometry. 

On the toric side, the sequence of papers \cite{A,BGS,CS,E,GH,GKS,LMP1,TW} led to proving that a closed Richardson variety $\Rvw$ is toric (with respect to a possibly bigger torus than $T$) if and only if $[v,w]$ is $S_3$-free (i.e., misses the $S_3$ Bruhat poset, also called a crown) if and only if $[v,w]$ 
is a lattice. In the simply laced cases $\Rvw$ is smooth toric if and only if $[v,w]$ is a Boolean lattice. Examples of $T$-stable $(\cp)^d$'s inside the flag variety were constructed in \cite{Ohn}.

On the Poisson side, the full flag variety $G/B_+$ admits the so-called standard Poisson structure, \cite{GY,L,STS}
\[
\pi := \sum_{\al \in \Delta_+} \chi(e_\alpha) \wedge \chi(f_\alpha),
\]
where $\Delta_+$ denotes the set of positive roots of $\g := \Lie(G)$, $\{e_\alpha \mid \alpha \in \Delta_+\}$ and $\{f_\alpha \mid \alpha \in \Delta_+\}$ are sets of dual root vectors with respect to the Killing form and $\chi : \g \to H^0( G/B_+, T (G/B_+) )$ is the infinitesimal $\g$-action on $G/B_+$. The $T$-action on $G/B_+$ keeps the Poisson structure $T$-invariant, and therefore transforms symplectic leaves to symplectic leaves. The $T$-orbits of symplectic leaves of $(G/B_+, \pi)$ are precisely the open Richardson varieties $\oRvw$, \cite[Example 4.9]{EL}, \cite[Theorem 0.4]{GY}. 

The Poisson degeneracy loci \cite{P} (subschemes where the Poisson structure is of rank $\leq 2k$) play an integral role in capturing the Poisson geometry of Fano manifolds through the Bondal conjecture asserting that those schemes have sufficiently high dimension \cite{Bondal,GP}. In  \cite{CMY} it was proved that the reduced part of the 0-degeneracy locus of $G/B_+$ (the subscheme where the Poisson structure vanishes) is equal to the union of Richardson varieties $\Rvw$ for all pairs of Weyl group elements $v \leq w \in W$ satisfying any of the following equivalent conditions:
\begin{enumerate}
    \item $\dim\Ker(vw^{-1}+1) = l(w) - l(v)$, where $l$ is the standard length function on $W$.
    \item $(v w^{-1})^2 =1$ and $l_\Delta(v w^{-1}) = l(w) - l(v)$ where $l_\Delta$ denotes the reflection length function on $W$.
    \item $v = s_{\gamma_1} \dots s_{\gamma_d} w$ for some positive roots $\gamma_1, \ldots, \gamma_d$ such that $\gamma_j \bot \gamma_k, \forall j \neq k$, where $d := l(w) - l(v)$.
    \item For one, and thus for each reduced word $(i_1, \ldots, i_{l(w)})$ of $w$, there exists a reduced subword $(i_1, \ldots , \widehat{i}_{p_1} , \ldots , \widehat{i}_{p_d} , \ldots , i_{l(w)})$ with value $v$ 
    such that $\beta_{p_{j}} \bot \beta_{p_{k}}$, $\forall j \neq k$ for the roots $\{\beta_1, \ldots, \beta_{l(w)} \} = \Delta_+ \cap w(- \Delta_+)$ given by $\beta_k: = s_{i_1} \ldots s_{i_{k-1}}(\al_{i_k})$ where $\al_i$ are the simple roots of $G$ and $s_i$ the corresponding reflections of $W$.
\end{enumerate}
Denote by $\GCR(W)$ the collection of pairs $(v,w) \in W^2$ with $v \leq w$ satisfying any of the above properties, where ECR stands for {\em{extended covering relation}}, 
see Sect. \ref{sec2} for details. It was proved in \cite{CMY} that if any of those conditions holds, then $\Rvw \simeq (\cp)^{l(w) - l(v)}$.

Comparing the toric and Poisson geometric treatments of $G/B_+$ gives that the reduced part of the 0-degeneracy locus of the flag variety $G/B_+$ is a union of special toric Richardson varieties $\Rvw$ and raises the question of how we can geometrically characterize those Richardson varieties $\Rvw$ that appear in the Poisson degeneracy locus of $G/B_+$. We settle this problem in the present note:

\begin{theorem}\label{main}
Let $v, w \in W$ and assume $v\leq w$. Then $$\Rvw \simeq (\cp)^{l(w)-l(v)} \iff (v,w) \in \GCR(W).$$
\end{theorem}

The backward implication was proved in \cite[Theorem B]{CMY}; the first main result in this note is the forward implication, which is harder. 

Up to now, symplectic foliations of Poisson varieties have been primarily studied from a set theoretic perspective \cite{EL,GY,LY} by describing their symplectic leaves individually as orbit closures. The upshot of Theorem \ref{main} is that we can now provide an algebro-geometric characterization. Firstly, recall that a complete Poisson subvariety of a Poisson variety $(X, \pi)$ is a subvariety that is a union of symplectic leaves, i.e., is stable under all Hamiltonian flows of $X$ (everywhere by a ``subvariety'' we mean a closed subvariety). By the above mentioned results of \cite{EL,GY}, the $T$-stable complete Poisson subvarieties of $(G/B_+, \pi)$ are the (arbitrary) unions of closed Richardson varieties $\Rvw$. The closure relation for Richardson varieties implies that the irreducible $T$-stable complete Poisson subvarieties of $(G/B_+, \pi)$ are precisely the closed Richardson varieties. Therefore, Theorem \ref{main} implies the following:

\newpage

\begin{corollary}\label{Pgeom}
\hfill
\begin{enumerate}
\item[(i)] The $T$-stable complete Poisson submanifolds of the standard Poisson structures on the full flag varieties $(G/B_+, \pi)$ that are isomorphic to powers of $\cp$ are precisely the closed Richardson varieties $\Rvw$ for $(v,w) \in \GCR(W)$. 
\item[(ii)] The reduced part of the 0-degeneracy locus of $(G/B_+, \pi)$ is the union of the $T$-stable complete Poisson subvarieties of $(G/B_+, \pi)$ that are isomorphic to powers of $\cp$. Its irreducible components are the maximal ones under inclusion.
\end{enumerate}
\end{corollary}

In the spirit of Bondal's conjecture \cite{Bondal}, we ask whether in general one can see the degeneracy loci of Poisson Fano varieties purely algebro-geometrically:

\begin{question}  Let $(X,\pi)$ be a Poisson Fano variety with an equivariant action of an algebraic group $T$. Is it true that the irreducible components of the reduced parts of the Poisson degeneracy loci of $(X,\pi)$ of level $k$ are precisely the maximal $T$-stable complete Poisson subvarieties of $(X,\pi)$ that belong to a given class $\mathcal{C}_k$ of projective varieties?
\end{question}

The class $\mathcal{C}_k$ should depend only on a family $\mathcal{F}$ to which the  Poisson Fano variety $(X,\pi)$ belongs. For instance in the case of $k=0$ and full flag varieties, $\mathcal{C}_0$ consists of powers of the projective line. 

The proof of Theorem \ref{main} is based on the following result on arbitrary $T$-stable $(\cp)^d$'s in $G/B_+$, which is the second main result of this note and is of independent interest:

\begin{theorem}\label{orthocell}
Let $X \subseteq G/B_+$ be a $T$-stable subvariety isomorphic to $(\cp)^d$. 
\begin{enumerate}
    \item[(i)] There exist $x \in W$ and pairwise orthogonal positive roots $\gamma_1,\ldots, \gamma_d$ such that 
$$X^T=\{ x \prod_{i \in I}s_{\gamma_i}B_+:I \subseteq [d] \},$$ where $[d]:=\{1,\ldots,d\}.$
    \item[(ii)] The variety $X$ is the closure of a single $T$-orbit. 
    \item[(iii)] If $K \subseteq T$ denotes the kernel of the action on $X$, then $T/K$ is a $d$-dimensional torus acting on $X$.
    \item[(iv)] Writing $p_I=x \prod_{i \in I}s_{\gamma_i}B_+$, the irreducible closed $T$-stable curves in $X$ are precisely the root curves joining $p_I$ and $p_J$ such that $|I \triangle J|=1$. 
    \item[(v)] We also have $$d \leq r_{\perp}(\Delta):=\max \{|\Omega|: \Omega \subseteq \Delta_+ \, \text{is pairwise orthogonal}\}.$$
\end{enumerate}
\end{theorem}

We note that $r_{\perp}(\Delta)$ is precisely the cardinality of the Kostant cascade of roots \cite{Kostant} of the Lie algebra of $G$. 

\subsection*{Acknowledgments} 
The authors thank the Institute of Mathematics and Informatics of the Bulgarian Academy of Sciences, where they first met during the Mathematics Days in Sofia conference.
The research of M.Y. was supported by NSF grant DMS-2601486, Bulgarian Science Fund grant KP-06-N92/5, the Ministry of Education and Science grant DO1-239/10.12.2024, and Simons Foundation grant SFI-MPS-T-Institutes-00007697. 

\subsection*{AI statement} The authors used large language models, including ChatGPT-6.1 and Claude Fable 5.2, as exploratory tools in connection with the problems studied in this paper. Their use included brainstorming, discussion of possible approaches, and preliminary checking of ideas. The authors take full responsibility for the content of the paper.

\sectionnew{Background and notation}
\label{sec2}
In this section we gather background material on Lie groups and Richardson varieties. Throughout the paper $G$ will denote a connected simply connected complex semisimple Lie group with Lie algebra $\g$ and Cartan matrix $C:=(c_{ij})_{i,j=1}^n$. Let $B_{\pm}$ be a pair of opposite Borel subgroups of $G$ and $T := B_{+}\cap B_{-}$ be the corresponding maximal torus of $G$. Denote $ \mathfrak{b}_{\pm}:= \Lie(B_\pm)$ and $ \mathfrak{t}:= \Lie(T)$. The corresponding sets of roots, positive roots, negative roots, and simple roots will be denoted by 
\[
\Delta, \quad \Delta_+, \quad \Delta_-, \quad \mbox{and} 
\quad \Pi := \{ \alpha_i \mid 1 \leq i \leq n \},
\]
respectively, thought of as subsets of the character lattice of $T$ and $\mathfrak{t}^{*}$. 
 Let $d_1, \ldots, d_n$ be the positive integers symmetrizing the Cartan matrix $C$ such that the subcollection of those corresponding to any simple factor of $G$ are relatively prime. Denote by $( \cdot , \cdot)$ the symmetric non-degenerate bilinear form on $\mathfrak{t}^{*}$ such that 
\[
(\alpha_i, \alpha_j) = d_i c_{ij}.
\] 

For any root $\beta$, denote the corresponding root subgroup $U_{\beta} \subset G$, as the unique closed connected $T$-stable subgroup of $G$ with Lie algebra $\mathfrak{g}_{\beta}$ (cf. \cite[Proposition 8.1.1]{Springer}).

Denote the Weyl group of $G$ by $W$. The reflection corresponding to a root $\beta$ will be denoted by $s_\beta$, and in the case of simple roots, by $s_i:=s_{\alpha_i}$.
For a spanning set $R \subseteq \Delta$ of $\mathfrak{t}^{*}$, the reflection length of an element $w \in W$ with respect to $R$, denoted $l_R(w)$, is the minimal $r \in \Nset$  such that $w$ can be written as a product of $r$ reflections $s_{\beta_1}, \ldots , s_{\beta_r}$ with $\beta_1, \ldots , \beta_r \in R$. 
When $R:= \Pi$, we recover the usual length function on $W$, $l_{\Pi}(w) = l(w)$. At the other extreme for $R:= \Delta$ (or what is the same $R:= \Delta_+$), we obtain the \emph{reflection length} of $w$ (also called absolute length or rank)
\[
l_\Delta(w).
\]
We will use the first terminology. 

The full flag variety $G/B_+$ is stratified into the \emph{open Richardson varieties} \eqref{oRvw}  
\[
G/B_{+} = \bigsqcup_{v \leq w} \oRvw .
\]
The \emph{closed Richardson variety} $\Rvw$ is the Zariski closure of $\oRvw$ in $G/B_+$:
\[
\Rvw = \bigsqcup_{v \leq v' \leq w' \leq w} \mathring{\mathcal{R}}_{v'}^{w'}.
\]
It is irreducible and has dimension $l(w)-l(v)$. 

\bth{ECM} \cite[Theorem A]{CMY} The following conditions for a pair of Weyl group elements $v \leq w \in W$ with $d := l(w) - l(v)$ 
are equivalent:
\begin{enumerate}
    \item $\dim\Ker(vw^{-1}+1) = l(w) - l(v)$.
    \item $(v w^{-1})^2 =1$ and $l_\Delta(v w^{-1}) = l(w) - l(v)$ where $l_\Delta$ denotes the reflection length on the Weyl group of $G$.
    \item $v = s_{\gamma_1} \dots s_{\gamma_d} w$ for some positive roots $\gamma_1, \ldots, \gamma_d \in \Delta_+$ such that $\gamma_j \bot \gamma_k, \forall j \neq k$.
    \item For one, and thus for each reduced word $(i_1, \ldots, i_{l(w)})$ of $w$, there exists a reduced subword $(i_1, \ldots , \widehat{i}_{p_1} , \ldots , \widehat{i}_{p_d} , \ldots , i_{l(w)})$ with value $v$ 
    such that $\beta_{p_{j}} \bot \beta_{p_{k}}$, $\forall j \neq k$ for the roots $\beta_k: = s_{i_1} \ldots s_{i_{k-1}}(\al_{i_k})$, $1 \leq k \leq l(w)$.
\end{enumerate}
\eth

Let $\GCR(W)$ be the collection of all pairs of Weyl group elements $v \leq w$ satisfying any of the above four equivalent conditions. This set is an extension of the covering relation $\lessdot$ of the Bruhat order on $W$ in the sense that 
\[
(w, w) \in \GCR(W) \quad \mbox{and} \quad (v,w) \in \GCR(W), \; \; \forall v \lessdot w.
\]

By \cite[Theorem B]{CMY} for all $(v,w) \in \GCR(W)$ the following hold:
\begin{enumerate}
    \item[(i)] The closed Richardson variety $\Rvw$ is isomorphic to a power of $\cp$
    \[
     \Rvw \simeq (\cp)^{l(w)-l(v)}.
     \] 
     \item[(ii)] The Bruhat interval $[v,w]$ is isomorphic to the power set of $\{1, \ldots, l(w) - l(v)\}$, i.e., is a Boolean lattice.
\end{enumerate}

\sectionnew{Proof of Theorem \ref{orthocell}}
\label{sec3}

Firstly, we recall that we have
\begin{equation}\label{t-fixed}
(G/B_+)^T=\{uB_+: u \in W\} 
\end{equation}
and that there is a bijection $$W \longrightarrow (G/B_+)^T, \, u \mapsto \dot{u}B_+,$$
where $\dot{u}$ is a representative of $u$ in the normalizer $N_G(T)$ of $T$. 

We will need some auxiliary results. A standard fact about torus actions on powers of projective lines is given by the following.

\begin{lemma}\label{torus-coordinates}
Let an algebraic torus $S$ act on $Y=(\cp)^d$. There exist $g_1,...,g_d \in PGL_2$ and characters $\chi_1,\ldots, \chi_d \in X^*(S)$ such that in the projective coordinates $z_i'=g_i(z_i)$, the action is given by 
$$t \cdot (z_1',...,z_d')=(\chi_1(t)z_1',\ldots, \chi_d(t)z_d'),$$ for every $t\in S$ and every $(z_1',...,z_d') \in Y$, where $c \cdot \infty = \infty$ for $c \in \mathbb{G}_m$. In particular $\mathbb{G}_m^d$ is $S$-stable and every point of $\{0,\infty\}^d$ is $S$-fixed. If every $\chi_i$ is non-trivial, then $Y^S=\{0,\infty\}^d$.

\end{lemma}

\begin{proof}

It suffices to show that the action is factorwise and that on each factor $S$ fixes two distinct points. Indeed, sending those points to $0$ and $\infty$ makes every element of $S$ act on that factor by multiplication by a scalar. 

Let $\pi_i: Y \longrightarrow \cp$ be the $i$-th projection. Choose a fiber  $F=\pi_i^{-1}(b)$ and a point $a\neq b$. The set $\{(t,z) \in S \times F:  \pi_i(t\cdot z)=a\}$  is closed. Its image $Z:=\{t \in S: a \in \pi_i(t \cdot F)\}$ under the proper projection $S \times F \longrightarrow S$ is therefore closed. 

Since $e \notin Z$, the complement $S\setminus Z$ is non-empty. It is also open and therefore dense since $S$ is irreducible. For $t \in S\setminus Z$ the map $z \mapsto \pi_i(t \cdot z)$ sends the connected proper variety $F$ into $\cp \setminus \{a\} \simeq \mathbb{A}^1$, so it is constant. Hence, for $z, z' \in F$, the two maps $t \mapsto \pi_i(t \cdot z)$ and $t \mapsto \pi_i(t \cdot z')$ from $S$ to $\cp$ agree on the dense subset $S \setminus Z$. Since $\cp$ is separated, the two maps agree on all of $S$. Since $i$ and $b$ were arbitrary, $S$ acts factorwise. 

It remains to prove that $S$ fixes two distinct points on each factor. If $\sigma_i:\cp\to Y$ is a section of $\pi_i$, the induced action on the $i$-th factor is $(t,z)\mapsto\pi_i(t\cdot\sigma_i(z))$. This map is algebraic, and the factorwise property makes it independent of the choice
of $\sigma_i$.

If this action is trivial, any two distinct points are fixed. Otherwise, its kernel is a proper closed subgroup of $S$, and since the torsion points of $S$ are Zariski dense, some $t_0\in S$ of finite order acts non-trivially. A non-identity M\"obius transformation of finite order is given by a diagonalizable non-scalar matrix, so $t_0$ has exactly two fixed points. Since $S$ is abelian, it permutes them, and since $S$ is connected, it fixes each of them.

Choose $g_i\in PGL_2$ sending the two fixed points to $0$ and $\infty$. In the coordinate $z_i'=g_i(z_i)$, each $t\in S$ acts by $z_i'\mapsto\chi_i(t)z_i'$ with $\chi_i(t)\in\Cset^\times$. Since $\chi_i(t)=t\cdot 1$, the map $\chi_i$ is regular, and the action law makes it multiplicative. Hence $\chi_i\in X^*(S)$. In homogeneous coordinates the action is $[u:v]\mapsto[\chi_i(t)u:v]$, so the formula holds on all of $Y$, with $c\cdot\infty=\infty$.

This formula shows that $\mathbb G_m^d$ is $S$-stable and that every point of $\{0,\infty\}^d$ is $S$-fixed. If every $\chi_i$ is non-trivial, then $0$ and $\infty$ are the only fixed points on each factor, and since $S$ acts factorwise, $Y^S=\{0,\infty\}^d$.

\end{proof}

\begin{remark}\label{Xisclosed}
Let $\iota:X \hookrightarrow G/B_+$ be the inclusion. Suppose $X \simeq (\cp)^d$. Then it is proper over $\Cset$. The graph of $\iota$ is closed in $X \times G/B_+$ because $G/B_+$ is separated. The projection $X \times G/B_+ \longrightarrow G/B_+$ is proper, so it sends the graph to a closed set. The image is $\iota(X)$. Therefore $X$ is closed in $G/B_+$. 
\end{remark}

We will also use the following fact.

\begin{lemma}\label{zero-coords}
Let $X \subseteq G/B_+$ be $T$-stable. Let $\varphi: X \rightarrow (\cp)^d$ be an isomorphism. For any $xB_+ \in X^T$, the coordinates in Lemma \ref{torus-coordinates} can be chosen so that $\varphi(xB_+)=(0,\ldots,0)$ and the characters $\chi_1,\ldots,\chi_d$ are distinct roots in $x(\Delta_-)$.  
\end{lemma}

\begin{proof}
Transport the action of $T$ through $\varphi$ and denote $p=\varphi(xB_+)$. Since $\varphi: X \rightarrow (\cp)^d$ is an isomorphism, the inclusion  $X \hookrightarrow G/B_+$ induces a $T$-equivariant injection of tangent spaces $$T_p((\cp)^d) \hookrightarrow T_{xB_+}(G/B_+).$$  
Choose coordinates as in Lemma \ref{torus-coordinates}, and write $z_i$ for the resulting coordinates. The weights of $T$ on $T_{xB_+}(G/B_+)$ are distinct roots in $x(\Delta_-)$. Since the $i$-th coordinate direction of $T_p((\cp)^d)$ has weight $\pm \chi_i$, and all the weights of $T_{xB_+}(G/B_+)$ are nonzero, each character $\chi_i$ is nontrivial. 

Since $xB_+ \in X^T$, the point $p$ is fixed by $T$, and Lemma \ref{torus-coordinates} gives $p \in \{0, \infty\}^d$. By replacing coordinates $z_i$ with $z_i^{-1}$, whenever $p_i=\infty$, we may assume that $p=(0,\ldots,0)$. 

Replacing a coordinate $z_i$ with $z_i^{-1}$ changes the corresponding character from $\chi_i$ to $-\chi_i$. After these changes, the characters $\chi_1,\ldots, \chi_d$ are the tangent weights at $(0,\ldots,0)$. The $T$-equivariant injection above shows that they are distinct roots in $x(\Delta_-)$, as desired.
\end{proof}

We omit the proof of the following lemma (cf. \cite[\S 8.2]{Graham-Kreiman}), which describes the $T$-stable curves in the full flag variety.

\begin{lemma}\label{root-curve}
For $y \in W$ and $\beta \in y(\Delta_-)$, the curve $\overline{U_{\beta}yB_+} \simeq \cp$ is the unique irreducible closed $T$-stable curve through $yB_+$ whose tangent weight there is $\beta$. Its fixed points are $yB_+$ and $s_{\beta}y B_+$.  Every irreducible closed $T$-stable curve through $yB_+$ arises in this way.
\end{lemma}

We call the curves $\overline{U_{\beta}yB_+}$ \emph{root curves}, joining $yB_+$ and $s_{\beta}y B_+$.

\subsection{Proof of Theorem \ref{orthocell} (i)}

Let $d=\dim X$ and fix an isomorphism $\varphi: X \longrightarrow (\cp)^d$. If $d=0$, the results of the theorem are immediate. We can therefore assume that $d > 0$. 

Transport the action of $T$ to $(\cp)^d$ through $\varphi$. Lemma \ref{torus-coordinates} guarantees that $X^T$ is non-empty since every vertex of $\{0,\infty\}^d$ is fixed. Let $xB_+ \in X^T$. Using Lemma \ref{zero-coords}, we may choose coordinates so that $\varphi (xB_+)=(0,\ldots,0)$ and the action of $T$ on all of $(\cp)^d$ is given by
\begin{equation}\label{coord-action}
a \cdot (z_1,\ldots, z_d)=(\chi_1(a)z_1,\ldots, \chi_d(a)z_d)    
\end{equation}
where $\chi_1, \ldots, \chi_d$ are distinct roots in $x(\Delta_-)$.

For $I \subseteq [d]$, let $q_I \in \{0, \infty\}^d$ denote the point whose $i$-th coordinate is $\infty$ when $i \in I$ and $0$ otherwise. Each character $\chi_i$ is non-trivial, so the induced $T$-action on the $i$-th copy of $\cp$ has precisely two fixed points, $0$ and $\infty$. Therefore $$X^T=\{\varphi^{-1}(q_I):I \subseteq [d] \}.$$

Now, since $\varphi^{-1}(q_I)$ is $T$-fixed, using (\ref{t-fixed}), we can write $\varphi^{-1}(q_I)=y_IB_+$, for $y_I \in W$. In particular we have $y_{\emptyset}=x$. 

We first show that swapping the coordinate $i$ between $0$ and $\infty$ changes the Weyl-group label by left multiplication by $s_{\chi_i}$, independently of the remaining coordinates. For ease of notation we write $r_i$ instead of $s_{\chi_i}$.

For $I \subseteq [d]$ and $i\in [d]$, consider the line $L_{I,i} \subseteq (\cp)^d$ consisting of the points $(z_1,...,z_d)$ such that $z_i$ varies freely and whenever $j\neq i$, then  $z_j=\infty$ for $j\in I$, and $z_j=0$ otherwise. The inverse image $\varphi^{-1}(L_{I,i})$ is an irreducible $T$-stable curve in $X$, and hence in $G/B_+$. Its two fixed points are $y_IB_+$ and $y_{I\triangle\{i\}}B_+$.

The curve $\varphi^{-1}(L_{I,i})$ is closed in $G/B_+$ since $L_{I,i}$ is closed in $(\cp)^d$ and $X$ is closed in $G/B_+$ (recall Remark \ref{Xisclosed}). Its tangent weight at $y_IB_+$ is $\chi_i$ if $i \notin I$, and $-\chi_i$ if $i \in I$.

By Lemma \ref{root-curve}, the fixed points of the curve $\varphi^{-1}(L_{I,i})$ are $y_I B_+$ and $r_iy_IB_+$. We therefore get 
\begin{equation}\label{coordinate-flip}
r_iy_I=y_{I\triangle\{i\}}, \, \text{for all } \, I \subseteq [d] \, \text{ and } \, i \in [d].    
\end{equation}

This yields $y_{\{i \}}=r_ix$, since $y_{\emptyset}=x$. Therefore, for any fixed indices $i \neq j$, we also get 
$$r_jr_ix=r_jy_{\{i\}}=y_{\{i,j\}}=r_iy_{\{j\}}=r_ir_jx.$$
We therefore have $r_ir_jx=r_jr_ix$, which then gives us $r_ir_j=r_jr_i$.

Since the reflections $r_i$ commute, iterating (\ref{coordinate-flip}) gives $$y_I=\left(\prod_{i \in I}r_i \right)x, \, (I \subseteq [d]).$$

Define $\gamma_i=-x^{-1}(\chi_i)\in \Delta_+$. The $\chi_i$ are distinct roots in $x(\Delta_-)$ and the root system is reduced, so they are pairwise nonproportional. Therefore, since $r_i$ commute, we get that $(\chi_i,\chi_j)=0$, for $i\neq j$. Since the Weyl group preserves the inner product, we get $(\gamma_i,\gamma_j)=0$. Therefore $\gamma_1,\ldots,\gamma_d$ are pairwise orthogonal. Moreover, $r_i=s_{\chi_i}=xs_{\gamma_i}x^{-1}$, thus 
$$y_I=x \prod_{i \in I}s_{\gamma_i}$$
and part (i) of Theorem \ref{orthocell} holds.

\subsection{Proof of Theorem \ref{orthocell} (ii).}

Keep the notation and coordinates from the proof of part (i). The characters $\chi_i=-x(\gamma_i)$ are non-zero and pairwise orthogonal for the form $(\cdot, \cdot)$, which is positive definite on the real span of $\Delta$. Hence they are linearly independent over $\Rset$, and in particular over $\Zset$, in $X^*(T)$. 

Consider the homomorphism $$\eta: T \longrightarrow \mathbb{G}_m^d, \, t \mapsto (\chi_1(t),\ldots,\chi_d(t)),$$
which induces a map on character lattices $\eta^*:{\Zset}^d \longrightarrow X^*(T)$, sending the $i$-th standard basis vector to $\chi_i$. The map $\eta^*$ is injective and therefore $\eta$ is surjective. Indeed, $\eta(T)$ is a closed subtorus in $\mathbb{G}_m^d$ and if $\eta(T) \neq \mathbb{G}_m^d$, then some nontrivial character of $\mathbb{G}_m^d$ would be trivial on $\eta(T)$, which would imply that $\ker(\eta^*)$ is non-zero, contradicting the injectivity of $\eta^*$. 

Consider the point $\mathfrak{1}:=\varphi^{-1}(1,\ldots,1) \in X$. By (\ref{coord-action}) we have $t \cdot (1,\ldots,1) = \eta(t)$, so the surjectivity of $\eta$ gives $$T \cdot \mathfrak{1}= \varphi^{-1}(\mathbb{G}_m^d).$$
Since $\mathbb{G}_m^d$ is open and dense in $(\cp)^d$, the orbit of $\mathfrak{1}$ is open and dense in $X$. Now, $X$ is closed in $G/B_+$ (recall Remark \ref{Xisclosed}), so the closure of $T \cdot \mathfrak{1}$ in $G/B_+$ equals $X$. This proves (ii).

\subsection{Proof of Theorem \ref{orthocell} (iii).}

Continuing with the notation from (ii), we have $$K=\ker(\eta) = \bigcap_{i=1}^{d} \ker (\chi_i).$$
Indeed, if $t \in \ker (\eta)$, then $\eta(t)=(1,\ldots,1)$ and (\ref{coord-action}) implies that $t$ acts trivially on $(\cp)^d$, and hence on $X$, which means that we have $t\in K$. Conversely, if $t\in K$, then $t\cdot \mathfrak{1} = \mathfrak{1}$, and therefore $\eta(t)=(1,\ldots,1)$.

Since $\eta$ is surjective, we get a bijective homomorphism from $T/K$ to $\mathbb{G}_m^d$, which is an isomorphism since we work in characteristic $0$. Thus $T/K$ is a torus of dimension $d$. The action of $T$ on $X$ factors through a faithful action of $T/K$, and under the identification $T/K \simeq \mathbb{G}_m^d$ (and the isomorphism $\varphi$) it becomes the coordinatewise action of $\mathbb{G}_m^d$ on $(\cp)^d$. This proves (iii). 

\begin{remark}
\begin{enumerate}
    \item[(a)] Note that parts (ii) and (iii) imply that $X$ is a toric variety for the torus $T/K$. More precisely, via $\varphi$ and the isomorphism $T/K \simeq \mathbb{G}_m^d$ induced by $\eta$, the variety $X$ is equivariantly isomorphic to $(\cp)^d$ with the coordinatewise $\mathbb{G}_m^d$-action. In particular, $X$ is a smooth projective toric variety for the torus $T/K$.
    \item[(b)] Note that $K$ need not be connected. For example, if $G=SL_2$ and $X=G/B_+ \simeq \cp$, then $K=\ker(\alpha)=\{\pm 1\}$, where $\alpha$ is the positive root. 
\end{enumerate}

\end{remark}

\subsection{Proof of Theorem \ref{orthocell} (iv).}

By the proof of part (iii), the $T$-orbits of $X$ correspond under $\varphi$ to products obtained by choosing $\{0\}$, $\{\infty\}$, or $\mathbb{G}_m$ in each factor. There are finitely many such orbits, so every irreducible closed $T$-stable curve is the closure of a one-dimensional orbit. These closures are precisely the coordinate lines $\varphi^{-1}(L_{I,i})$.  

By the proof of part (i), the $T$-fixed points of $\varphi^{-1}(L_{I,i})$ are $p_I$ and $p_{I \triangle \{i\}}$, and Lemma \ref{root-curve} identifies $\varphi^{-1}(L_{I,i})$ as a root curve. Conversely, every pair $p_I, p_J$ with $|I\triangle J|=1$ occurs in this way. Finally, the root curve joining such a pair is unique. Indeed, if two root curves through $yB_+$ have the same other endpoint, then $s_{\beta}y= s_{\gamma}y$ for $\beta, \gamma \in y(\Delta_-)$. Thus $s_{\beta}= s_{\gamma}$, and $\beta=\gamma$ because the root system is reduced and $y(\Delta_-)$ does not contain pairs of opposite roots. Hence every root curve joining $p_I$ and $p_J$ with $|I\triangle J|=1$ is one of the lines $\varphi^{-1}(L_{I,i})$ and in particular lies in $X$. This proves (iv).

\subsection{Proof of Theorem \ref{orthocell} (v).}

This follows directly from part (i) of the present theorem and the definition of $r_{\perp}(\Delta)$.

\sectionnew{Proof of Theorem \ref{main}}
\label{sec4}

Let $v, w \in W$ with $v\leq w$, and set $d=l(w)-l(v)$.  If $(v,w) \in \GCR(W)$, then  \cite[Theorem B]{CMY} proves that $\Rvw \simeq (\cp)^d$. For the converse, assume $X:=\Rvw \simeq(\cp)^d$. By Theorem \ref{ECM}(3), it suffices to write 
\begin{equation}\label{vSOSw}
v=s_{\gamma_1}\cdots s_{\gamma_d}w,
\end{equation}
for some positive pairwise orthogonal roots  $\gamma_1, \ldots, \gamma_d \in \Delta_+$.

If $d=0$ then $v=w$ and the result is immediate. We therefore assume that $d > 0$. 

We first recall the fixed-point description (cf. \cite[Lemma 1 (2)]{BrionLakshmibai}) 
\begin{equation}\label{richardson-fixed}
X^T=\{uB_+: v \leq u \leq w\}.    
\end{equation}

Note that for $u \in W$ and $\gamma \in \Delta_+$ we have  

\begin{equation}\label{bruhat-root-sign}
u < u s_{\gamma} \Longleftrightarrow l(u) < l(u s_{\gamma}) \Longleftrightarrow u(\gamma) \in \Delta_+.
\end{equation}

Applying Theorem \ref{orthocell}, part (i), to $X$, we can find an element $x \in W$ and pairwise orthogonal positive roots $\gamma_1,\ldots,\gamma_d,$ such that 
$$X^T=\{ x \prod_{i \in I}s_{\gamma_i}B_+:I \subseteq [d] \}.$$

We will show that we can choose our basepoint to be $v$. But first write $A_I:=\prod_{i \in I}s_{\gamma_i}$. Orthogonality implies the reflections $s_{\gamma_i}$ commute, so $A_IA_J=A_{I \triangle J}$ and $A_I^2=1$.

By (\ref{richardson-fixed}), $vB_+\in X$, so $vB_+=xA_JB_+$ for some subset $J \subseteq [d]$. Distinct Weyl group elements label distinct fixed points, so $x=vA_J$. Using the properties of $A_J$ above, we have $$xA_I=vA_JA_I=vA_{J\triangle I}.$$

Since $I \mapsto J\triangle I$ permutes the subsets of $[d]$, we get
$$X^T=\{ v \prod_{i \in I}s_{\gamma_i}B_+:I \subseteq [d] \}.$$

Writing $x_I=v \prod_{i \in I}s_{\gamma_i}$, and applying (\ref{richardson-fixed}), we see that the fixed points of $X$ are indexed by the Bruhat interval
$[v,w]=\{x_I: I \subseteq [d]\}$.

Since $vs_{\gamma_i}\in [v,w]$ and $vs_{\gamma_i} \neq v$, we have $v < vs_{\gamma_i}$. Therefore, using (\ref{bruhat-root-sign}), we get $v(\gamma_i) \in \Delta_+$.

For $i\notin I$, orthogonality implies $A_I(   \gamma_i)=\gamma_i$  and therefore $x_I(\gamma_i)=v(\gamma_i)\in \Delta_+$, therefore using (\ref{bruhat-root-sign}) again, we have 
\begin{equation}\label{addone}
x_I<x_Is_{\gamma_i}=x_{I\cup \{i\}}.    
\end{equation}

Write $w=x_{J'}$ for some ${J'} \subseteq [d]$. If ${J'}\neq [d]$, choose $i\notin {J'}$. Then applying the inequality (\ref{addone})  to our $x_{J'}$ gives $w=x_{J'}<x_{{J'}\cup \{i\}}$, which contradicts the maximality of $w$. Therefore $J'= [d]$ and
\begin{equation}\label{wmax}
w=v\prod_{i=1}^{d}s_{\gamma_i}.
\end{equation}

Now set $\delta_i:=v(\gamma_i) \in \Delta_+$ and note that these are pairwise orthogonal roots (since $\gamma_i$ are such). Also, we have $s_{\delta_i}=vs_{\gamma_i}v^{-1}$, and by multiplying these relations we therefore get $$s_{\delta_1}\cdots s_{\delta_d}=vAv^{-1},$$ where we write $A=s_{\gamma_1}s_{\gamma_2}\cdots s_{\gamma_d}$.

We now have
$$s_{\delta_1}\cdots s_{\delta_d}w = (vAv^{-1})(vA)=vA^2=v,$$ where we used the fact that $A^2=1$, since the reflections $s_{\gamma_i}$ commute.

This proves (\ref{vSOSw}) as desired and hence the theorem.

\end{document}